\documentclass{article}

\usepackage[T1]{fontenc}
\usepackage[utf8]{inputenc}

\usepackage[leqno]{amsmath}
\usepackage{amssymb}
\usepackage{mathtools}
\usepackage{latexsym}
\usepackage{bussproofs}
\usepackage{accents}
\usepackage{graphicx}
\usepackage{natbib}

\usepackage{url}
\DeclareUnicodeCharacter{03C9}{\ensuremath{\omega}}

\DeclareUrlCommand\email{\urlstyle{rm}}

\vfuzz \hfuzz
\usepackage[vskip=\baselineskip]{quoting}

\renewenvironment{quote}[1]%
{\begin{quoting}#1}
{\end{quoting}}

{\begin{quoting}#1}
{\end{quoting}}

\usepackage[protrusion=true,stretch=10,shrink=10]{microtype}

\usepackage{etoolbox}

\usepackage{geometry}%
\newcommand*\GetTextWidth[3][\normalfont\normalsize]{{#1%
    \settowidth{#2}{abcdefghijklmnopqrstuvwxyz}%
    \setlength{#2}{0.03193#2}%
    \addtolength{#2}{0.44961pt}%
    \setlength{#2}{#3#2}%
    \global#2=#2}}

\newlength\bringhurstwdt
\GetTextWidth{\bringhurstwdt}{72}
\usepackage[german,dutch,french,english]{babel}
\useshorthands{"}

\newcommand{\dutch}[1]{\foreignlanguage{dutch}{#1}}

\newcommand{\french}[1]{\foreignlanguage{french}{#1}}
\newcommand{\german}[1]{\foreignlanguage{german}{#1}}

\newcommand{\weg}[1]{}

\newcommand{\rotiota}[2]{\rotatebox[origin=c]{180}{$#1\iota$}}
\newcommand{\rotatediota}{{\mathpalette\rotiota\relax}}

\newcommand{\formatabr}[1]{{\MakeUppercase{#1}}}

\newcommand{\cnrs}{\formatabr{cnrs}}

\newcommand{\ens}{\formatabr{ens}}
\newcommand{\iv}{\formatabr{iv}}
\newcommand{\snd}{\formatabr{snd}}

\makeatletter
\g@addto@macro \normalsize {%
\setlength{\abovedisplayskip}{\baselineskip}%
\setlength{\belowdisplayskip}{\baselineskip}%
\setlength{\abovedisplayshortskip}{\baselineskip}%
\setlength{\belowdisplayshortskip}{\baselineskip}%
\setlength{\jot}{0pt}
}
\makeatother

\DeclareOldFontCommand{\bf}{\normalfont\bfseries}{\mathbf}

\title{Predicativity and parametric polymorphism of Brouwerian implication}
\author{Mark van Atten\thanks{\snd{} (\cnrs{}\,/\,Paris \iv{}),
1 rue Victor Cousin,
75005 Paris,
France. 
As of September 1,
2018:
Archives Husserl (\cnrs/\ens{}),
45 rue d'Ulm,
75005 Paris,
France.
Email:
\mbox{vanattenmark@gmail.com}}}

\begin{document}


\maketitle

\begin{abstract}
\noindent A common objection to the definition of intuitionistic implication
in the Proof Interpretation
is that it is impredicative.
I discuss the history of that objection,
argue that in Brouwer's writings predicativity of implication is ensured through
parametric polymorphism of functions on species,
and compare this construal with the alternative approaches 
to predicative implication
of
Goodman,
Dummett,
Prawitz,
and
Martin-Löf.
\end{abstract}

\section{Impredicativity and intuitionistic implication}

A definition is impredicative
if it defines a member of a totality
in terms of that totality itself.
There is a familiar contrast between two
kinds of impredicative definitions,
as illustrated by these two examples:
\begin{enumerate}
\item \(k = \text{the largest number in } \{1,2,3\}\);
\(k\) is a member of that totality.

\item Let 
\(s\)
be an infinite set.
According to the power set axiom,
there exists a set 
\(\mathcal{P}(s)=\rotatediota x(\forall y(y \in x \leftrightarrow y \subseteq s))\).
Here 
\(\forall y\)
is a quantification over the totality of all sets;
\(\mathcal{P}(s)\)
is a member of that totality.
\end{enumerate}

In the first,
the totality in question
is enumerated
and this is one way of defining
each of its elements
prior to defining 
\(k\);
the impredicative definition of 
\(k\) 
serves to single out
a certain element of that totality,
but is not essential to the introduction of that element.
On the other hand,
for an infinite set
\(s\),
there is no alternative
to characterising
the set of all its subsets
by relating it to the totality of all sets,
because
there is no way to generate
all subsets of a given infinite set
from below.%
\footnote{That is,
inductively.
It is possible to define the 
\textit{species}
of subsets of a given infinite set,
but not to generate its extension;
see p.~\pageref{L005} below on Brouwer’s acceptance of the species of all species of real numbers.
In his criticisms of Cantorian set theory in his dissertation,
Brouwer rejects arbitrary exponentiation,
but does not go into the notion of power set.
See on this point 
\citealt[p.113]{Dalen1999}. 
See also 
\citealt{Poincare1909}
(thanks to Thierry Coquand for the suggestion).}

For constructivists this means that
the first definition is acceptable as a definition
of a mathematical object,
whereas
the second
is not.%
\footnote{The Brouwer-Kripke Schema can be used to construct an intuitionistic analogue of sorts to \(\mathcal{P}(N)\)
\citep{Dalen1977}.
For an analysis of that schema,
see
\citealt{Atten-forthcomingD}.}
More generally,
if at a given moment,
the only definition available of a certain object
is an impredicative one that refers
to a totality for which we at that moment have no construction,
then this definition cannot be used to guide a construction process of that object.
To construct the defined object,
we would first have to construct the totality that its definition refers to,
but as this totality contains
the object we are in the process of constructing,
we find that we can only complete that process if we already have completed it.
In such cases the circularity of the impredicative
definition is vicious.


Of Brouwerian intuitionism
it has been argued 
that various of its definitions 
are impredicative in the vicious sense,
and hence not constructive:
the propositional connective of implication,
the principle of induction,
certain sequences in the so-called Theory of the Creating Subject (Troelstra’s Paradox),
the theory of ordinals,
certain species in Brouwer.
I will here focus on 
the clause for implication in the Proof Interpretation;
it will turn out that this requires also a discussion of species,
and sheds some light on the theory of ordinals.%
\footnote{For a discussion of the principle of induction, 
see \citealt{Atten-forthcomingC}; 
for Troelstra’s Paradox, 
\citealt{Atten-forthcomingB}.}
An often cited formulation of it is that
in
\textit{Constructivism in Mathematics}
by Troelstra and van Dalen:
\begin{quote}
A proof of 
\(A \rightarrow B\) 
is a construction which permits us to transform any proof of 
\(A\)  
into a proof of 
\(B\).
[Troelstra and van Dalen 1988, vol.1, p.9]
\end{quote}

In this formulation,
one will want to understand ‘construction’ 
as a
function 
\(f\),
because then the explanation what a proof of an implication is
will 
have the definiteness,
stability,
and
uniformity
characteristic of a function:
for each argument there will be a determinate value,
to a given argument at any time the same value will be assigned,
and to all arguments a value is assigned in the same way.

How is the domain of such a function
\(f\) 
to be understood?

It cannot be given
by 
an inductive definition of all proofs of 
\(A\),
as there can be no such definition.
In particular,
proofs of 
\(A\)
may themselves contain
the implication 
\(A \rightarrow B\),
for example in this way:
\begin{prooftree}
\AxiomC{$(A \rightarrow B) \rightarrow A$} 
\AxiomC{$A \rightarrow B$}
\RightLabel{$\rightarrow$-E}
\BinaryInfC{$A$}
\end{prooftree}
Thus a prior explanation
of
\(A \rightarrow B\)
would be required,
rendering
the definition circular.

But if the domain is to
be
all proofs of
\(A\)
in an absolute sense that can not be given by an inductive definition,
then,
so the claim goes,
the definition of such a function 
\(f\)
will take a form
that renders it
impredicative:
\begin{quote}
\(f\)
is a function such that, 
for any \(x\)
in 
the totality of all intuitionistic proofs,
if
\(x\) is a proof of 
\(A\),
then 
\(f(x)\) is a proof of 
\(B\).
\end{quote}
Any specific definition of this form
will define an individual proof 
of 
\(A \rightarrow B\)
by referring to a totality to which
it belongs,
and thus be impredicative.

Intuitionists 
consider the notion of proof to be open-ended,
not only epistemically 
(at no moment do we know all possible proofs)
but ontologically,
and 
hence they deny that there is such a thing
as the totality of all intuitionistic proofs
\citep[pp.148-149]{Brouwer1907}.
There is only a 
\textit{growing} 
universe of mathematical objects and proofs.
For those who share the criticism that intuitionistic implication is impredicative,
this intuitionistic conception of the universe would of course be part of what makes
the whole enterprise incoherent:
Intuitionists cannot have both a universe that is growing and a notion of implication
that demands that it is a totality.
Hence the intuitionistic denial that the universe is a totality
only makes the question
more urgent
whether there
are understandings of the notions
of proof,
function,
and 
domain,
such that
they allow for an intuitionistically coherent
reading of the clause for implication in the Proof Interpretation.

It might seem possible to circumvent the problem
by defining
\(A \rightarrow B\)
to mean
‘From the assumption that \(A\) is true,
a proof of \(B\)
can be obtained’,
the idea being that by making a mere assumption,
talk of proofs is avoided altogether.%
\footnote{This is done in ‘Basic Logic', e.g., \citealt{Ruitenburg1993}.}
But on Brouwer’s notion of truth,
according to which
‘truth is only in reality
i.e.~in the present and  past experiences of consciousness’
\citep[p.1243]{Brouwer1949C},
the fundamental meaning of
‘Proposition 
\(A\)
is true’
is that a suitable mental construction 
that is correctly described by
\(A\)
has been carried out.
Brouwer allows the idealisation
that,
whenever the Subject has not actually carried out a given construction method,
but could do so in principle,
the construction it leads to may be counted among those that have been carried out.
But in each case an appeal to this must be made explicitly,
for otherwise there is no conceptual reduction
to what has been experienced.
The latter
remains the primary notion.
Either way,
there is no such thing as a mere assumption.
To assume that \(A\) is true is to assume that a construction for  \(A\) has been carried out
(perhaps in an idealised sense).

A kind of 
‘weak impredicativity’\label{L001}
may also be
associated with the Proof Interpretation.
The term is taken from a draft note of 1970 by Gödel,
not on the Proof Interpretation,
but on his computable functionals of finite type:
\begin{quote}
In particular, there exist functions of lower type which, 
within $T$,
can only be defined in terms of functions of higher types. 
This is a
kind
of impredicativity. True, 
it is only one of those weak
impredicativities
that are admitted even in Princ[ipia] Math[ematica] 2nd ed.~p.[XL],%
\footnote{[The typescript erroneously has `XI', as had the original
publication of Gödel's Russell paper 
\citep[p.134]{Godel1944}.]}
ff.
In our proofs of the axioms of $T$ this impredicativity appears in the
fact
that the concept of reductive proof may itself occur in reductive
proofs 
(just as in Heyting's logic the general concept of proof may
occur in a proof).
\citep[p.218]{Atten2015}
\end{quote}
In the case of intuitionistic implication,
this would pertain to the composition of the function 
\(f\),
which 
may depend on  proofs of implications  up to a definite,
but arbitrarily high,
type,
as
a proof of a statement may refer to proofs of more complex statements
\citep[p.233]{Troelstra1990}.
To take a simple example,
\(f\) may be 
\( h\circ g\),
with
\(g\) 
raising the type and 
\(h\) 
appropriately lowering it again.
It may seem then as if
the existence of functions of a lower type
presupposes
the existence of functions of higher types.
I will return to this problem
after the discussion of
strong impredicativity.

\section{History of the objection}

It has been claimed 
\citep[section 2]{Dean.Kurokawa2015}
that 
the thought that intuitionistic implication is impredicative
goes back to Gödel’s
lecture 
titled 
‘The present situation in the foundations of mathematics’,
presented
at a meeting of the Mathematical Association of America
in Cambridge,
Massachusetts on December 30, 1933
\citep{Godel1933o}.
Contrary to the organisers’ plan,
Gödel left the text unpublished,
but it was included in the
\textit{Collected Works} 
in 1995.
And the suggestion that the objection of impredicativity
is made there seems natural;
in that text Gödel questions the constructive character
of the intuitionistic explanation what a proof of a
negation consists in,
which is a particular case of that for implication,
and his scepticism turns on the observation
that
this clause involves quantification over all possible proofs.
In this section,
I will first argue that,
all the same,
the objection that Gödel is actually making in 1933 is a different one,
and then conjecture that it was shortly after  Yale lecture of 1941
that the thought  
that intuitionistic implication is impredicative
must have occurred to him.
I will then consider the question who was the first to make that objection in print.

Gödel’s 1933 lecture is concerned with the question
of a constructive consistency proof for
classical arithmetic.
In considering what should count
as constructive mathematics,
Gödel there argues against accepting impredicative definitions,
and insists on inductive definitions.
Gödel discusses the prospects
for a consistency proof for classical arithmetic
using intuitionistic logic,
then best known from Heyting's formalisation
‘\german{Die formalen Regeln der intuitionistischen Logik}’
\citep{Heyting1930a,Heyting1930b,Heyting1930c},
as well as Heyting's Königsberg lecture
of 1930,
‘\german{Die intuitionistische Grundlegung der Mathematik}’,
published as
\citealt{Heyting1931}.
When writing the text for his talk, 
Gödel was able
moreover
to consult Heyting’s manuscript on intuitionism meant for their joint
book, 
eventually published by Heyting alone in 
\textit{\german{Mathematische Grundlagenforschung}}
\citep{Heyting1934}.
Heyting had sent the manuscript of his part to Gödel in August 1932 
\citep[p.54]{Godel2003a},
more than a year before Gödel's talk.

In
‘\german{Die formalen Regeln der intuitionistischen Logik}’,
Heyting says about 
\(a \supset b\)
that it means
‘If \(a\) is correct [‘richtig’], then so is \(b\)’
\citep[p.5]{Heyting1930a}.
In 
his Königsberg lecture of 1930,
Heyting did not discuss implication in general,
but about the special case of negation he said:
‘The proposition 
“\(C\) is not rational”
means the expectation
that from the assumption that
\(C\)
is rational,
a contradiction can be derived’
\citep[p.113]{Heyting1931}.
Gödel knew that lecture well,
because he had been in the audience,
and in 1932
he reviewed the published version for the
\textit{\german{Zentralblatt}}
\citep{Godel1932f}.
But it is in 
\textit{\german{Mathematische Grundlagenforschung}}
that 
one finds the first
formulation of the general clause for implication.%
\footnote{In Brouwer's work one finds various applications of it
\textit{\french{avant la lettre}};
see 
\citealt[section 3.1]{Atten2009b}.}
On p.14, 
Heyting describes the construction required to prove an implication
$a \supset b$
as one  
‘which from each proof of $a$ leads to a proof of $b$’.%
\footnote{‘\german{die aus jedem Beweis für a zu einem Beweis für b führt}’.}
This is a formulation that,
unlike the earlier ones,
makes it explicit that
this is 
an operation 
on proofs of $a$ and not a proof from 
the assumption $a$.
I do not know
whether Heyting's book of 1934 contained significant changes compared
to the manuscript he had sent Gödel in 1932.
But if so, it will not have been on this particular point,
as Heyting’s correspondence with fellow intuitionist Freudenthal 
shows that Heyting had
obtained the reading of implication 
as
‘I possess a construction that derives from every proof of 
\(a\) 
a proof of 
\(b\)’
shortly after the Königsberg lecture
\citep[pp.206–207]{Troelstra1983a}.
It is this understanding that is assumed in Gödel's criticism
in the 1933 lecture.

The principles in Heyting's formalisation that have
Gödel's special interest 
are those for  
‘absurdity’,
that is,
intuitionistic negation.
But Gödel goes on to argue that this notion is
not constructive in his sense,
and hence of no use for a 
constructive consistency proof of classical arithmetic.
The problem he sees 
is that their intuitionistic explanation
involve a reference to the totality of all constructive proofs.
The example he gives is
\[
p \supset \neg\neg p
\]
which,
he says, 
means %
‘If $p$ has been proved,
then the assumption $\neg p$
leads to a contradiction’.%
\footnote{Gödel says this example shows 
‘the character of the axioms assumed by Heyting’,
but in Heyting’s paper this one is not an axiom but a
theorem 
\citep[p.49, 4.3]{Heyting1930a}.
In a footnote,
Heyting says that he 
\emph{had} 
considered it an axiom,
until Glivenko pointed out to him that it was derivable.}
Gödel says
that
these axioms are not about constructions on a substrate of numbers
but rather on a substrate of proofs,
and therefore the example may be explicated as
‘Given \textit{any}
proof for a proposition 
$p$,
you can construct a reductio ad absurdum for the proposition
$\neg p$’. 
He then comments that
\begin{quote}
Heyting’s axioms concerning absurdity and similar notions
[…]
violate the principle,
which I stated before,
that the word
‘any’ can be applied only to those totalities for which we
have
a finite procedure for generating all their elements […]
The
totality of all possible proofs certainly does not possess
this
character,
and nevertheless the word ‘any’ is applied to this
totality in Heyting’s axioms 
[…] 
Totalities whose
elements
cannot be generated by a well-defined procedure are in some
sense
vague and indefinite as to their borders.
And this objection
applies particularly to the totality of intuitionistic proofs
because of the vagueness of the notion of constructivity.
Therefore this foundation of classical arithmetic by means of the
notion of absurdity is of doubtful value.
\citep[p.53]{Godel1933o}
\end{quote}
A draft of this passage in Gödel’s archive
does not quite end with a rejection
of Heyting’s logic.
Instead,
it reflects:
\begin{quote}
Therefore you may be doubtfull [sic] as to the correctness of the
notion of
absurdity and as to the value of a proof for freedom from
contradiction by means of this notion.
But nevertheless it may be
granted that this foundation is at least more satisfactory
than the
ordinary platonistic interpretation […]%
\footnote{\citealt[Series IV, box 7b, folder 26, document \#040113, 
page 22]{Godel.Papers}.
Work of Kurt Gödel used with permission. 
Unpublished Copyright (1934-1978) 
Institute for Advanced Study. 
All rights reserved by Institute for Advanced Study.}
\end{quote}
Either way,
the doubt about,
or objection to,
the notion of absurdity
immediately generalises to implication as such.

%

It is
remarkable,
given the construction of Gödel's talk,
in which  the discussion
of the intuitionistic logical connectives
is preceded by an argument
against the use impredicative definitions for foundational
purposes,
that the objection Gödel puts forward
is not that Heyting's principles for absurdity
are impredicative,
but that they are vague.
Impredicativity of course entails
constructive undefinability and in that sense
vagueness,
and it is possible that
Gödel had seen the problem of impredicativity
but thought that,
in the context of a consistency proof
that is looked for because of its epistemic interest,
vagueness is the more important thing to note,
even if impredicativity is the cause of it.

But vagueness may arise on other grounds,
and indeed 
to me it seems that the vagueness Gödel finds 
problematic in 1933 is not that caused by impredicativity, 
but
by open-endedness. 
After all,
when Gödel writes that
\begin{quote}
For the totality of all
possible proofs certainly does not possess this character
[of being an infinity that is generated by a finite procedure], 
and nevertheless
the word ‘any’ is applied to this totality in Heyting's axioms
\end{quote}
the second clause 
is evidently not presented as an explanation of the first.
But if one’s objection were that of impredicativity,
one would do just that.
One would argue that 
to accept the application of the word
‘any’  to the totality  of all
possible proofs,
an application that occurs in Heyting’s explanation
of what a proof of a negation consists in,
is to accept the existence of a
proof that is  impredicatively defined,
and therefore of a proof that will never
be generated by a finite procedure;
and that,
by implication,
the totality of all
possible proofs does not possess the character
of being an infinity that is generated by a finite procedure,
and hence is not constructive in the sense Gödel requires.
But this is not an argument one finds in his lecture.
In fact,
while Gödel does point out the problem,
from a constructive point of view, 
with
impredicatively defined integers, 
he does not say anything at all about
impredicatively defined proofs.
The first incompleteness theorem,
on the other hand,
was  brought up 
earlier in the lecture, 
and explained as meaning that
no single consistent formal system can embrace all formal proofs.
There is,
then, 
no finite generating procedure to exhaust all possible proofs.
And that the term
‘any’
as it figures in Heyting’s axioms
indeed must be taken to apply
to the totality of all possible proofs,
and cannot be explained via the notion of provability
in any single formal system,
was a result,
based on the second incompleteness theorem, 
that Gödel had 
contributed to Menger's
seminar in 1931/1932,%
\footnote{It was published in \citealt{Godel1933f}.}
well before writing his Cambridge lecture.%
\footnote{When Gödel in the lecture goes on to say that
‘this objection applies particularly to the totality of intuitionistic proofs
because of the vagueness of the notion of constructivity’ 
\citep[p.53]{Godel1933o},
he is pointing to a reinforcement of the problem,
not the problem as he sees it come about.}


Also
in his subsequent discussions in 
the ‘Lecture at Zilsel’s’ 
of 1938
and 
the lecture at Yale
‘In what sense is intuitionistic logic constructive?’
of 1941,
Gödel emphasises that intuitionistic provability 
cannot mean provability in a particular formal system,
and is therefore not constructive in what he there takes to be the strictest sense,
but he does not thematise impredicativity of any logical operation
(\citealt[pp.100-102]{Godel1938a}; 
\citealt[p.190]{Godel1941}).

As long as Gödel held on to the 
‘concrete’ 
notion of constructivity
outlined in each of the lectures of 1933, 1938, and 1941,
a notion founded on the idea of finite generating procedures, 
his argument that
intuitionistic logic is vague on account of the incompleteness theorem
must have convinced him that it is of no epistemic value in consistency proofs. 
But it no longer could when, 
in the year after the
Yale lecture, 
Gödel came to accept,
like the intuitionists, 
abstract notions as part of his concept of
constructivity;
the suggestion occurs already in his Princeton
lectures on intuitionism,
which ended on May 1,
whereas the Yale lecture had been delivered on April 15
\citep[section 11.3.5.2]{Atten2015}.
In Gödel’s case,
the abstract notion accepted as constructive was
not the intuitionistic notion of proof,
but that of computable functional of finite type,
which he showed can replace the former 
when limited to its application in arithmetic. 
(Gödel describes the necessity of admitting abstract notion
in print only in his
\textit{Dialectica}
paper of 
\citeyear{Godel1958}; 
in the introduction of that paper Gödel is,
in effect, 
also
addressing his former self.)
On the other hand,
Gödel continued to hold,
unlike the intuitionists,
that
there exists
a totality of all proofs.%
\footnote{In the 1972 version of the 
\textit{Dialectica} 
paper,
Gödel introduced the notion of 
‘reductive proof',
by which he hoped to 
arrive at an
interpretation of intuitionistic logic
that
‘in no way presupposes Heyting's and that, moreover, 
it is constructive and
evident in a higher degree than Heyting's. 
For it is exactly the elimination of such 
\textit{vast generalities} 
as 
“any proof” 
which makes for greater evidence and constructivity’
\citep[p.276n, emphasis mine]{Godel1972};
see also 
\citealt[sections 11.3.5.7 and 11.3.5.9]{Atten2015}.}
But while the predicativity of his old notion
of constructivity was ensured by its definition,
this was not  so for the widened one.
This became 
(and remained) 
a matter of great concern to Gödel.
If the
thought that implication in Heyting's Proof Interpretation
might be impredicative had not
occurred to him before, 
it will have then.


Not much time later,
between November 1942 and May 1943,
Gödel wrote his paper 
‘Russell's mathematical logic’,
published in 1944.
There,
Gödel 
comes to speak of the ‘self-reflexivity’ of
impredicatively defined properties 
\citep[p.139]{Godel1944}%
\footnote{Russell, 
in his paper 
‘On some difficulties in the theory of transfinite numbers and order
types’
\citeyearpar{Russell1906},
which Gödel refers to in his own paper,
had used the simpler term
‘reflexive’.}
and points out that
this self-reflexivity
need not be problematic in general,
but,
when it occurs in definitions
that are meant to be constructive,
renders them viciously circular.
However, 
he does not add that,
specifically,
the property of being an intuitionistic proof of an implication 
is problematic.

The first manuscript of Gödel’s
in which it is stated that intuitionistic logic
is impredicative
is an intermediate version of around 1970 of the revised 
\textit{Dialectica}
paper,
and
it is also included in the last,
1972 version,
posthumously published in 
\textit{Collected Works}.%
\footnote{Gödel Papers, 9b/144, item 040454, pages 13, end of note (j),
9b/144, item 040456, pages 6, end of note (j),
and item 040459, pages 6, end of note (j);
and
\citealt[p.274 note f]{Godel1990}.}
By then,
it had of course been stated in print elsewhere.

Although Gödel,
in print,
remained silent on the question,
his use of the term
‘self-reflexivity’
in relation to impredicativity in the Russell paper
seems to have struck a chord.
When 
Heyting published his book
\textit{Intuitionism. An Introduction}
\citeyearpar{Heyting1956},
probably the most influential publication
on intuitionism ever,
it was immediately reviewed
by Sigekatu Kuroda,
in the
\textit{Journal of Symbolic Logic}.
That review may well be the first
publication of,
in effect, 
the objection that intuitionistic implication
is impredicative:
\begin{quote}
The proof of the fan theorem in this book is strikingly
brief,
although it is essentially
the same as Brouwer’s original proof of 1926
(\textit{Mathematische Annalen} vol.97, p. 66, Th. 2).%
\footnote{[\citealt{Brouwer1927B}]}
By means of the dialogue method the author discusses a
certain peculiarity
of the proof of the fan theorem,
namely that it depends on inferences,
i.e.,
constructions,
which are assumed to have been executed previously.
The comment of the author
(or a person in the dialogue)
about it is as follows:
‘If we are well aware that the
hypothesis of a theorem consists always in the assumption of
a previous execution of
some construction,
we can offer no objection against the use of considerations
about
the way in which such a construction can be performed as a
means of proof.’ 
The
reviewer fully agrees to this comment.
However,
the assumed construction may contain
a construction by the application of other fan theorems,
in which further previous
constructions are assumed,
and so on.
In this way the assumption of a fan theorem
may contain,
so to speak,
in a nested manner many other fan theorems.
The same
situation arises when we wish to clarify the intuitionistic
concept of negation
(cf.~7.1.1),
and of implication,
as well as that of a function.
The self-reflexive character
of these could never be completely clarified by any formal
method.
[…]
Nevertheless,
some clear analysis
of the proof of the fan theorem would be necessary,
because no proof may be self-reflexive.
\citep[p.369]{Kuroda1956}
\end{quote}
It is not wholly clear to me to what extent
Kuroda takes an intuitionistic
implication
to be different from
a theorem based on a hypothesis;
on Brouwer's conception of truth,
however,
they are the same,
as to assume that 
\(P\)
is to assume that
\(P\)
has been demonstrated.
The identification of assuming the antecedent of an implication
and
assuming that a prior construction has been effected
is also clearly shown by Heyting's clause for implication in
the book under review:
\begin{quote}
The implication 
\(p \rightarrow q\) 
can be asserted, 
if and only if we possess a construction 
\(r\), 
which, 
joined to any construction proving 
\(p\) 
(supposing that the latter be effected), 
would automatically effect a construction proving 
\(q\). 
\citep[p.98]{Heyting1956}
\end{quote}
Perhaps the reason why Heyting 
in chapter 3 does not present the fan theorem in logical
notation
was that he was there focused on its
mathematical content,
whereas logic is discussed only in the penultimate chapter of his book;
the same remark can in fact be made about all theorems outside the chapter on logic.

Be that as it may,
Kuroda signals
the same problem for both.
Although he does not invoke impredicativity here,
it clearly is the problem he is describing,
using instead the same term that Gödel had associated with
impredicative definitions and its problem for constructivists
in his Russell paper,
‘self-reflexivity’.

Kuroda's review came out in 1956,
after he had
passed the academic year 1955-1956 at the 
Institute for Advanced Study,
where he had had discussions with Gödel and also with Kreisel
\citep[p.248]{Kuroda1958}.%
\footnote{‘The author has begun this work while at the Institute for Advanced
Study in Princeton, to which the author would like to express his gratitude
for its giving him an opportunity to ponder on the foundations of mathematics,
specifically to Professors Kurt Gödel and Georg Kreisel for many
strict criticisms and valuable discussions on this work at that time.’
\citep[p.2]{Leopoldt1975}
}
Kuroda had published four paper intuitionistic logic 
and its philosophical aspects 
before
\citeyearpar{Kuroda1947,Kuroda1948,Kuroda1949,Kuroda1951}.
One readily imagines therefore that the question of the meaning of intuitionistic logic was among the topics
of his conversations with Gödel and Kreisel;
however,
that has so far not been documented.

Indeed, 
it is Kreisel who soon after Kuroda’s
review
does use the term ‘impredicative’
to characterise the intuitionistic notion of proof,
in his 1958 review of \textit{Wittgenstein’s Remarks on
the Foundations of Mathematics}:
\begin{quote}
\mbox{}[Intuitionism] goes beyond finitism because it makes
statements concerning all possible constructions,
which certainly do not constitute a concrete totality.
[…]
undecided propositions,
and even implications between such propositions,
may be used as premisses in implications,
i.e.~one makes assertions which involve an hypothetical proof,
namely a proof of the premise,
though the totality of all proofs is not concretely specified
[…] Finitist mathematics does not use the general notion of a
constructive proof at all,
in fact it might be said to avoid
logical inferences
(\textit{which involve an impredicative concept of proof})
because it is restricted to purely combinatorial operations.
\citep[pp.147-148; emphasis mine]{Kreisel1958a}
\end{quote}

Later Kreisel added a twist: 
in 
‘Mathematical logic’
\citep{Kreisel1965} 
he accepts a certain theory as constructive,
even though on the interpretation he proposes it is not predicative,
on the ground that the Proof Interpretation is not predicative either:
\begin{quote}
We give here a formulation for 
intuitionistic logic
which is proof theoretically equivalent to
\(R\) 
[Feferman's 1964 system of predicative analysis].
[…] [N]ote that,
for the interpretation of Section 2
[i.e., the Proof Interpretation],
this theory 
does
\emph{not}
have predicative character at all,
because the intuitionistic logical operators
are themselves defined self-reflexively.
\citep[p.176, 3.531]{Kreisel1965}
\end{quote}
(Note the use of ‘self-reflexive’.)
Myhill
in ‘The formalization of intuitionism’
of 1968
comments 
that it is an ingenious ad-hominem
argument (p.337),
and indeed it is not an argument that
settles the matter.

\section{Brouwerian implication is predicative and parametrically polymorphic}
 
Myhill goes on to say that
‘Brouwer’s own practice leaves the issue
[of predicative vs impredicative intuitionism]
quite undecided’.
That impression is based
on cases in which Brouwer seems
to countenance certain
impredicative definitions,
of which we will see examples below.
But
I will argue that 
this is only seemingly so,
and that Brouwer
did not accept impredicative definitions.
In particular,
he had notions
of 
proof,
function,
and
domain
that
allow for an understanding of implication
that is predicative.
Readers familiar with the literature on Martin-Löf's Constructive Type Theory
will see a similarity of the account presented in this section
with the conception of type theory as a predicative open-ended system
with an infinite series of expanding universes;
in particular,
with
Michael Rathjen's point  
‘that all functions that deserve to be
called effective must at least be definable in a way that is
persistent with expansions of the universe of types’ 
\citep[p.427]{Rathjen2009}.
My aim here is not to analyse such similarities
but to present,
for its own sake, 
the older Brouwerian wherewithal
that happens to bring this similarity about.

Brouwer held that
mathematics consists first of all in mental acts of construction
\citep{Brouwer1907},
and 
Heyting explained the relation between such constructions and proofs
as follows:
\begin{quote}
If mathematics consists of mental constructions,
then every mathematical theorem is
the expression of a result of a successful construction.
The proof of the theorem
consists in this construction itself,
and the steps of the proof are the same as the
steps of the mathematical construction.
\citep[p.107]{Heyting1958A}
\end{quote}
After certain acts of mental construction have been carried out,
these may
in reflection,
be objectified.
For Brouwer this is one form 
(the nonlinguistic one) 
of ‘second-order mathematics’:
the mathematical viewing of mathematical acts
\citep[pp.98,119n,173]{Brouwer1907}.
There is,
then, 
a primary sense of
proof as process,
and a secondary one
of proof 
as objectified process,
and the object constructed in the process  is
the whole of the constructed objects and relations between them
that make the proved proposition true;
a state of affairs.

In order to keep proofs in Brouwer’s sense,
that is
proofs in the sense of episodes in the mental life of the Creating Subject,
terminologically distinct
from proofs in the sense of linguistic or abstract objects outside
the mind,
we may propose to call Brouwer’s proofs
\textit{demonstrations}
\citep{Sundholm.Atten2008};
the ambiguity between demonstration as act
and
demonstration as object
is accepted and dealt with by adding the appropriate qualification.

The most memorable example of Brouwer’s use of 
demonstration objects is
in
his  paper
‘\german{Über Definitionsbereiche von Funktionen}’
\citep{Brouwer1927B},
in which he gives a demonstration of the Bar Theorem
and 
in 
footnote 8
adds this elucidation:
\begin{quote}
Just as,
in general,
well-ordered species are produced by means of the two
generating operations from primitive species
(see my paper [\citealt[p.~451]{Brouwer1927A}]),
so,
in particular,
mathematical demonstrations
[‘\german{Beweisführungen}’] 
are
produced by means of the two generating operations from null
elements and elementary inferences 
[‘\german{Elementarschlüssen}’]
that are immediately given in intuition
(albeit subject to the restriction that there always occurs
a last elementary inference).
These \emph{mental} 
mathematical proofs 
[‘\german{\textit{gedanklichen} Beweisführungen}’]
that in general contain
infinitely
many parts
[‘\german{Glieder}’]
must not be confused with their linguistic
accompaniments,
which are finite and
necessarily inadequate,
hence do not belong to mathematics.

The preceding remark contains my main argument against the
claims of 
Hilbert’s metamathematics\label{L002}. 
\citep[p.64n8]{Brouwer1927B}
\end{quote}

Brouwer says here that a 
\german{Beweisführung}
is a 
\textit{species},
one of Brouwer's two constructive analogues to a classical set,%
\footnote{The other one is the \textit{spread}, which in fact is a special case.}
and a well-ordered one at that.
What will be of particular interest for the present
purpose is not that individual demonstrations are
species,
but that,
being mathematical objects,
demonstrations may be collected
in a second-order species.
Brouwer defined species as 
\begin{quote}
properties supposable for mathematical entities previously acquired,
and satisfying the condition that,
if they hold for a certain mathematical entity,
they also hold for all mathematical entities which have been defined to be equal to it,
relations of equality having to be symmetric,
reflexive and transitive;
mathematical entities previously acquired for which the property holds
are called the elements of the species.
\citep[p.142]{Brouwer1952B}
\end{quote}
Here,
an object qualifies as  ‘previously acquired’ 
if either it has actually been constructed,
or a method has been indicated for constructing it.%
\footnote{%
\label{L004}
This explicit condition was absent from the first definition of species
(1918):
\begin{quote}
Unter einer Species erster Ordnung verstehen wir
eine Eigenschaft,
welche nur eine mathematische Entität besitzen kann,
in welchem Falle sie ein Element der Species erster Ordnung genannt wird.
\mbox{}\citep[pp.3-4]{Brouwer1918B}
\end{quote}
\parindent 15pt
It does not occur in the similar definition of 1925 either:
\begin{quote}
Unter einer Spezies erster Ordnung verstehen wir
eine 
(begrifflich fertig definierte)
Eigenschaft,
welche nur eine mathematische Entität besitzen kann,
in welchem Falle sie ein Element der Species erster Ordnung genannt wird.
\citep[p.245]{Brouwer1925A}
\end{quote}
\parindent 15pt
It seems to make its first appearance in handwritten additions to a manuscript related to the Berlin lectures, which where held in 1927:
\begin{quote}
\mbox{}[…]
moreover the admission
at each stage of this construction of mathematics
of properties which
can be presupposed for 
mathematical 
previously acquired mathematical  thought-entities,
as new mathematical thought-entities under the name of
species 
\citep[pp.483-484, trl.~modified]{Stigt1990}
[\dutch{‘…
en daarbij de toelating,
in ieder stadium van dezen opbouw der wiskunde,
van eigenschappen,
die voor reedsch verkregen mathematische denkbaarheden ondersteld kunnen worden,
als nieuwe mathematische denkbaarheden onder den naam van \emph{soorten}.’} 
\citep[p.482]{Stigt1990}]
\end{quote}
At the time,
Brouwer had the project of publishing the
Berlin lectures as a book
\citep[p.551]{Dalen2005},
but this addition may have been made at any time;
it was incorporated almost verbatim
in  
\citealt[p.23]{Brouwer1947}.}
Higher-order species 
were present from the beginning:
A species of order
\(n\)
has as elements 
mathematical entities or species of order
\(n-1\) 
\citep[p.4]{Brouwer1918B}.


An example of objects that count as all acquired
even though there are infinitely many of them
is of course the natural numbers,
given that for them
we possess a construction method.
Under this definition,
then,
the existence of an object precedes the existence of any species that it may become an element of,
and this rules out essentially impredicative definitions of elements of a species.
Curiously,
given his early interest in Poincaré and Russell,
Brouwer never pauses to state this;
Heyting later does
(\citealt[p.111]{Heyting1931}, 
\citealt[p.25]{Heyting1934}).%
\footnote{These passages confirm
Beeson's supposition 
\citep[p.52]{Beeson1985}
that
Heyting’s intention was to insist that species be defined predicatively.}
As we will see,
Brouwer  does make the remark that species, 
thus defined, 
cannot be elements of themselves.

But the definition of a species is well-motivated independently of its constructively 
welcome consequence that the definition of an element of a species is always predicative.
This independent motivation can be seen by considering that
the existence of a species 
\(S\) 
amounts to the existence of
an assignment,
to each existing object 
\(a\),
of
the hypothesis that it has the property 
\(P\).
As Brouwer says,
a species is 
a ‘property
\textit{supposable} 
for mathematical entities previously acquired’
(emphasis mine).
To accentuate the role of hypotheses
in Brouwer's species
is one of the virtues of Van Stigt’s
discussion of that notion
\citep[section 6.3.6]{Stigt1990};
I should like to add
that this assignment
is a
function.
Its domain is 
\textit{every} 
existing object
(all objects ‘previously acquired’) 
because for each of them the proposition that it has that property is meaningful
(whether that proposition turns out to be true, false, or undecided).

What the domain of the function can not include is
objects beyond those that have been previously acquired,
as otherwise
to construct the domain we should be able to effect a construction
out of non-acquired objects,
which is impossible;%
\footnote{This means also that the function cannot be a member of its own domain;
see the postscript to
\citealt{Atten-forthcomingB}.}
and the associated hypotheses would not be meaningful.
It would seem,
then,
that while this conception 
obviously allows for the existence of the species of the natural numbers,
there 
could be no such thing as, 
for example,
the species of real numbers.
There are more than finitely many of them, 
but no construction method for generating all,
and hence at no point they will all have been acquired.
But Brouwer adds the following comment:
\begin{quote}
With regard to this definition of species we have to remark
firstly
that,
during the development of intuitionistic mathematics,
some species will have to be considered as being re-defined time and again in the same way,
secondly that a species can very well be an element of another species,
but never an element of itself.
\citep[p.142]{Brouwer1952B}
\end{quote}
With every further acquisition of objects,
the function implicit in the definition of a species
that generates the hypotheses
\(P(a)\)
should be defined anew,
because this acquisition means that we want to extend
the domain of that function;
and the definition of the domain is a constitutive part
of the definition of a function.%
\footnote{So the need of redefinition in the case of certain species was implicit in the definition
of species as soon as it started to include the qualification ‘previously acquired entities’;
see footnote~\ref{L004}.}
But the way in which the new objects belong to the species should be the same as
that for the objects it already contains.
In the case of the species of real numbers,
each object that belongs to it does so because it has the property of being
a convergent choice sequence of rational intervals.
This warrants speaking of the
redefined species being the same species
as the previous one:
it is a growing object.

This redefinition of species along the different stages
of the Creating Subject’s activity is of course analogous to the situation in 
the later developed and more familiar Kripke models,
where
a predicate may be valuated differently
along a path through possible worlds
(‘evidential situations’),
with later valuations including the earlier ones.
By this I do not mean to make the farther-reaching suggestion
that Kripke models explain the intended meaning of intuitionistic logic.%
\footnote{Indeed it was not designed to meet that desideratum:
\begin{quote}
I do not think of ‘possible worlds’ as providing a
\emph{reductive}
analysis in any philosophically significant sense,
that is,
as uncovering the ultimate nature,
from either an epistemological or a metaphysical point of view,
of modal operators,
propositions,
etc.,
or as 
‘explicating them’. […] 
The main and the original motivation for the
‘possible worlds analysis’ 
–~and the way it clarified modal logic~–
was that it enabled modal logic to be treated by the same set theoretic techniques
of model theory that proved so successful when applied to extensional logic.
It is also useful in making certain concepts clear.
\citep[p.19n18]{Kripke1980}
\end{quote}
}
It is rather a stepping stone towards the idea that,
to make the notion of a species
such as Brouwer uses it
fully explicit,
one may associate to it 
an infinite sequence
that starts with a species in the strict sense
(i.e., with a fixed domain),
and in which
every further species is
obtained in a separate defining act of the Creating Subject,
as an extension of the previous one.

A special case of growing species
are those
for which
we have a
a method to generate denumerably many elements of it,
but also a method to extend any species
thereby obtained.
This is very close to
Russell’s 
notion of
‘self-reproductive processes and classes’.
\begin{quote}
\mbox{}[A]ccording to current logical assumptions,
there are what we may call self-reproductive processes and classes. 
That is, 
there are some properties such that, 
given any class of terms all having
such a property, 
we can always define a new term also having the property
in question. 
Hence we can never collect all the terms having the said
property into a whole […]
\citep[p.36]{Russell1906}
\end{quote}
This is from Russell's paper 
‘On some difficulties in the theory of transfinite numbers and order
types’ of 1906.
Brouwer's notebooks towards his dissertation,
and that dissertation itself,
defended in February 1907,
contain no indication that he was aware of Russell’s paper.
Be that as it may,
the
‘denumerably unfinished sets’
that he introduces in the dissertation
are very similar:
\begin{quote}
\mbox{}[H]ere we call a set
denumerably unfinished
if it has the following properties:
we can never construct in a well-defined way
more than a denumerable subset of it,
but when we have constructed such a subset,
we can immediately deduce from it,
following some previously defined mathematical process,
new elements which are counted to the original set.
But from a strictly mathematical point of view this set does not exist as a whole,
nor does its power exist;
however we can introduce these words here as an expression
for a known intention.%
\footnote{[‘Intention’ translates the Dutch ‘\dutch{bedoeling}’.]}
(\citealt[p.148]{Brouwer1907}, trl.~\citealt[p.82]{Brouwer1975})
\end{quote}
Brouwer's definition is motivated
by a preceding discussion why
Cantor's second number class
is not,
as a whole, 
a constructible entity;
it is a denumerably unfinished set.
Following this definition
are
the examples of
the totality of definable points on the continuum
and a fortiori the totality of all possible mathematical systems.
Brouwer's notebooks towards his dissertation
show that he made his first remark on
denumerably unfinished sets early on
(probably 1905),
while reading 
chapter 55, 
on projective geometry, 
in
Russell’s
\textit{Principles of Mathematics}
\citeyearpar{Russell1903}:%
\begin{quote}
One will never be able to resolve the whole mystery
of space and surfaces in such a way,
that there is
no hocus-pocus
to it anymore;
for the number of possible
buildings [i.e., constructions]
is denumerably unfinished,
hence not surveyable.
\citep[Notebooks 1904–1907, vol.II, p.32]{Brouwer19041907}%
\footnote{‘\dutch{Men 
zou het heele mysterie van ruimte en vlakken nooit kunnen ophelderen
zóó, 
dat er geen hokuspokus meer aan is; 
immers het aantal mogelijke gebouwen
is aftelbaar onaf, 
dus niet te overzien.’}}
\end{quote}
Later,
Dummett would discuss collections of that type,
with reference to Russell but not to Brouwer, 
as
‘indefinitely extensible concepts’.%
\footnote{First in ‘The philosophical significance of Gödel's Theorem’ 
\citep{Dummett1963},
without reference to Russell.
For a more elaborate discussion,
with the reference,
see   
\citealt[pp.316–319]{Dummett1991a}.
Various connections can be made also to the discussion
of such concepts
by Wright and Shapiro \citeyearpar{Shapiro.Wright2006}
in their 
‘All things indefinitely extensible’.}
It has been noticed that the concept of a denumerably unfinished set
more or less disappears from Brouwer’s writings after the introduction
of choice sequences 
\citep[p.115]{Dalen1999};
with the new concept of redefined species,
Brouwer could moreover subsume that of denumerably unfinished set under the latter.

Illustrative in this respect
is 
Brouwer’s assertion,
 in a paper of 1927,
of  the existence of
‘the species \(O\) of ordinal numbers’
\citep[p.487]{Brouwer1927A}.
Firstly,
because 
strictly speaking there is no such thing as
\textit{the}
species of ordinal numbers;
whenever we have defined an enumeration of ordinals,
this can be used to construct a new ordinal.
But 
it is the implicit redefinitions of the species that allow Brouwer 
to speak this way.%
\footnote{Thus,
I do not share the view expressed in
\citealt[p.127]{Gielen.Swart.Veldman1981}
that Brouwer's accepting this species \(O\) is an indication that
over the years he had become ‘more liberal’.}
Secondly, 
because he is not calling
\(O\)
a denumerably unfinished set here.
Earlier,
in a published reply to Mannoury’s review
of his dissertation,
Brouwer had called
the set of 
definable points on the continuum
denumerably unfinished;
and although
that is not the same set,
Brouwer here is
otherwise
repeating his reply to Mannoury.%
\footnote{See \citealt{Mannoury1908} and \citealt{Brouwer1908B}
(pp.162 and 170 of the reprints in \citealt{Dalen2001a}),
and, 
for further discussion of 
‘denumerably unfinished sets’,
\citealt[section 7.6.2]{Kuiper2004}.}

Similarly,
in the
\textit{Cambridge Lectures}
of 1946–1951
he can
accept
the species of all species of real numbers\label{L005} 
\citep[p.26n]{Brouwer1981A},%
\footnote{In contrast,
its editor,
Dirk van Dalen,
in a footnote \citep[p.26n18]{Brouwer1981A}  
comments that
‘Brouwer literally quotes the power set of the reals’.
(He continues by pointing out that Brouwer does not actually
need the power set for the example he gives on that page.)}
and he could,
a fortiori, 
have accepted
the species of all species of natural numbers.%
\footnote{Heyting \citeyearpar[p.195]{Heyting1962}
considered it
‘doubtful’
that these formed a species,
and preferred not to use the notion.}

Implicit redefinition of a species is also essential to considerations
whether some 
\textit{hypothetical}
object 
would belong to a certain species:
This is possible because,
on Brouwer's conception of existence, 
to reason
about a hypothetical object 
is to reason
under the hypothesis that that object 
has been acquired,
and hence that the species has been redefined accordingly.

And just as in the case of the species,
the function that is implicit in its definition
may have to be redefined time and again,
so will
every function defined on such a species.
The reason is
again that
the change in the species which is the domain
of the function
necessitates a redefinition of the function,
even though  the way the new function acts on the elements of its domain
is the same as for the old function.
These redefinitions of a function on a species
can be made explicit
as an infinite sequence of functions
parallel to the infinite sequence of species.

The one thing that needs to be verified
is that 
from the mere schemata of the species and the function defined on it
a sufficient justification can be found
for the constructibility of 
such functions
on 
\textit{arbitrary} 
elements of the species.
In the case of demonstration,
the most obvious and most frequently used property
is that a demonstration of 
\(A\) 
has conclusion 
\(A\).
If that is all the information that the function will use,
then it will indeed work for arbitrary demonstrations of
\(A\).
But in the case of the Bar Theorem,
Brouwer considered that for
\(A\) = ‘There is a (decidable) bar in the tree’ 
a further structural property
is known,
namely that every demonstration object of 
\(A\) 
can be put into a specific canonical form.

In Brouwer’s view,
as we saw,
demonstration objects of propositions 
are themselves objects of mathematics
and hence
may
become elements of species.
The species of demonstration objects of 
\(A\)
is not an example Brouwer gives explicitly.
But when we understand implication
as
\begin{quote}
There exists a function
\(f\) 
that transforms
any element of the species of demonstration objects of 
\(A\)
into an element of the species of demonstration objects of 
\(B\)
\end{quote}
\noindent
our understanding should reflect the fact
that both this function and
the species that is its domain
(and,
perhaps,
the species that is its range)
are 
\textit{growing} 
objects.
As constructions,
these take the form of infinite sequences. 
Hence,
that formulation  is to be taken as shorthand for
\begin{quote}
There exists an infinite sequence of
species of demonstration objects of 
\(A\) 
indexed by their stage of definition 
\(m\),
and an infinite sequence of 
functions
\(f_m\),
each extended by the next,
such that 
\(f_m\) 
transforms
any element
of the species of demonstration objects of 
\(A\) 
at stage 
\(m\) 
into
an element of the species
of demonstration objects of 
\(B\)
at stage 
\(m\).
{\renewcommand{\arraystretch}{1.2}
\[
f : A \rightarrow B \quad = \quad 
\begin{array}{cccccccc}
A_1 & \subseteq & A_2 & \subseteq & A_3 & \cdots & A_m & \cdots\\
\Bigg\downarrow\vcenter{%
\rlap{$f_1$}} 
&  & 
\Bigg\downarrow\vcenter{%
\rlap{$f_2$}}   
& & 
\Bigg\downarrow\vcenter{%
\rlap{$f_3$}}   
  & & 
\Bigg\downarrow\vcenter{%
\rlap{$f_m$}}   
 & \\
B_1 & \subseteq & B_2 & \subseteq & B_3 & \cdots & B_m & \cdots
\end{array}\]
}
\end{quote}

By thus taking ‘the species of all demonstration objects’,
and hence any subspecies 
‘demonstration objects of 
\(A\)’
to be a 
\textit{\french{façon de parler}} 
for an infinite sequence of ever larger species,
and a function on such a species to be
a 
\textit{\french{façon de parler}}  
for an infinite sequence of functions that
all have been, 
or will be, 
defined according to the same schema,
an understanding of implication becomes possible
that does not require
quantification over the totality of demonstration objects.%
\footnote{Independently,
a formal semantics showing
a structurally similar approach to intuitionistic implication
has simultaneously been developed in
\citealt{Tabatabai}.}

From the perspective of
modern type theory,
and the theory of programming languages,
one would not call it a 
\textit{\french{façon de parler}},
but rather an 
implementation,
namely,
of 
parametric polymorphism.%
\footnote{The notion was introduced in
a lecture course by
Christopher Strachey in 1967,
published as \citealt{Strachey2000}.
See also
\citealt{Cardelli.Wegner1985}.
Simultaneously with the first preprint of the present paper,
the preprint \citealt{Pistone} appeared,
which draws attention to the epistemological significance
of parametricity for the study of impredicativity,
there in the context of consistency arguments for impredicative systems.}
A function is
parametrically polymorphic if
its arguments need not each be of one fixed type,
but may come from a family of types,
on which the function however acts uniformly;%
\footnote{In intuitionistic mathematics,
that should not mean that the function cannot
inspect its arguments at all;
in the Bar Theorem,
where non-canonical proofs of the antecedent are to
be transformed into canonical form, 
on which then a proof of the consequent is based,
it is even essential that one can.
What cannot be done is to treat,
upon inspection,
different arguments instantiating the same polymorphic type differently.}
the types of the arguments are considered to be
parameters,
and with different instantiations of these parameters
the function takes a different shape.
A standard example is a function that takes
finite sequences of elements of arbitrary type
and returns the length of the sequence.
In the case of Brouwerian implication,
as reconstructed here,
a function 
\(f\)
that proves an implication
\(A \rightarrow B\)
is likewise 
parametrically polymorphic.
The family of argument types it can act on
consists of the species of 
demonstration objects of 
\(A\) 
at stage 
\(m\),
for all
\(m\);
the function 
\(f\)
has,
so to speak,
infinitely many types at once.
Depending on
\(A\),
the arguments of the function 
\(f\)
may themselves be polymorphic functions
(iteration of implication),
but all polymorphic functions used here are bounded
in that their domains are 
limited to demonstration objects 
that have been
‘previously acquired’.
This induces a stratification.
In particular,
then,
none of these polymorphic functions can be an element of its own domain.
The form of parametric polymorphism described here is
more expressive than unstratified predicative polymorphism
(where polymorphic type variables are always instantiated with non-polymorphic types),
and less expressive than impredicative polymorphism
(where polymorphic type  variables may be instantiated with any type).%
\footnote{It would be interesting to know more,
historically,
about Brouwer’s knowledge of Russell’s ramified theory of types.}
This is known as
‘finitely stratified polymorphism’ \citep{Leivant1991};
the finiteness comes naturally in Brouwerian intuitionism,
where mathematical construction acts are taken to proceed in
an 
\(ω\)-order,
so that any stage of the Creating Subject’s
activity
has only finitely many predecessors,
and hence by stage 
$n$
there exist only
the
finitely many species of 
demonstration objects of 
\(A\) 
at stage 
\(m<n\).

From this point of view,
the upshot of the present account,
then,
is that not only is intuitionistic logic an appropriate logic to reason about parametricity,%
\footnote{For an overview, see \citealt{Meseguer1989}. 
Of particular interest is also 
\citealt{McCarty1991},
which relates polymorphism to Brouwer’s notion of apartness.} 
but,
on a Brouwerian conception of logic at least,
it does itself arise from parametricity.
There is no circularity because
for Brouwer
logic depends on mathematics,
not the other way around,
and
parametricity 
is a mathematical phenomenon.%
\footnote{Conversely,
the notion of parametricity is broadly applicable in mathematics.
See~\citealt{HermidaReddyRobinson2014}.}

Heyting’s formalisation of species in 1930,
which predates the discussion of impredicativity in the Proof Interpretation,
is set up in a way
that,
in hindsight, 
has the effect of
specifically blocking
parametric polymorphism:
\begin{quote}
With this notation it is not possible
to extend the domain of a mapping
without changing the name of the mapping.
If one has defined
e.g.~the function
\(\sin\) 
on the domain of the real numbers,
and wishes to extend it to the domain of the complex numbers,
one from now on needs to
refer to that
which so far was called
\(\sin\)
by
\mbox{\(\sin \upharpoonright \textrm{(real numbers)}\)}
\citep[p.12, trl.~mine]{Heyting1930b}.%
\footnote{‘\german{Bei dieser Schreibweise ist es nicht möglich,
das Gebiet einer Abbildung zu erweitern,
ohne die Benennung der Abbildung zu ändern.
Wenn man z.B.~die Funktion
\(\sin\) 
für das Gebiet der reellen Zahlen definiert hat,
und man will sie auf das Gebiet der komplexen Zahlen erweitern,
so ist man genötigt,
dasjenige,
das bisher 
\(\sin\) 
hieß,
weiterhin durch 
\mbox{\(\sin \upharpoonright \textrm{(reellen Zahlen)}\)}
anzudeuten.’}}
\end{quote}

And Heyting in 1956 
defines species in such a way that
the need to redefine them,
together with functions defined on them,
is lost sight of,
or rather,
is suggested not to exist:
\begin{quote}
Definition 1. A species is a property which mathematical entities can be supposed to possess (L. E. J. Brouwer 1918, p. 4; 1924, p. 245; 1952, p. 142).

Definition 2. After a species S has been defined, 
any mathematical entity which has been 
\textit{or might have been} 
defined before S and which satisfies the condition S, 
is a member of the species S.
\citep[section 3.2.1, emphasis mine]{Heyting1956}
\end{quote} 
This is of importance for a historical understanding of the discussion
of impredicativity in intuitionism,
because at least until the publication
of \citealt{Brouwer1975}, 
Heyting’s book was a, 
if not the, 
major reference in the discussion.

The present analysis of Brouwer's species and functions on them
also allows for
the weak impredicativity that Gödel indicated
(see above, p.~\pageref{L001}).
Gödel was referring to the fact that 
the definition of a
computable functional
from a low type to another low type may 
refer to functionals of arbitrarily higher types;
for an implication 
this corresponds to the idea
that a proof of it may contain
proofs of statements that are arbitrarily more complex
than the antecedent or the consequent.
In his Russell paper,
Gödel had argued that,
if one accepts Russell and Whitehead’s axiom
that
functions occur in propositions only extensionally
(introduced  in the second edition of 
\textit{Principia}, 
pages xl-xli),
then definitions of functions exhibiting
this kind of impredicativity
‘are quite unobjectionable
even from the constructive standpoint […] provided that
quantifiers are always restricted to definite orders’
\citep[p.134]{Godel1944}.
On the other hand,
intuitionistically that axiom
would be unnatural.
A construction that depends on a function 
\(f\),
such as a construction that is a truth-maker of a proposition referring to
\(f\),
should be allowed to depend on 
how
\(f\)
is given to us.
It is,
however,
possible to define a function at a higher type
before certain functions at a lower type
that are defined so as to depend on the former,
as long as the higher-type function is
appropriately generic;
that is,
as long as it is guaranteed that the higher-type
function will also work as expected on
arguments that are constructed after it.

\section{Comparison with other predicative accounts of implication}

Other constructive accounts of implication
that,
like Brouwer’s,
employ the strategy 
of stratification
to give an account of implication
that does not quantify over a totality of all possible proofs,
are
the
‘theory of constructions’
of Kreisel 
\citeyearpar{Kreisel1962,Kreisel1965}  and Goodman 
\citeyearpar{Goodman1970,Goodman1973a,Goodman1973b}%
\footnote{For a recent discussion of this theory
(these theories), 
and further references,
see \citealt{Dean.Kurokawa2015}.}
and the
meaning-theoretical approaches of
Dummett \citeyearpar{Dummett1975}, 
Martin-Löf \citeyearpar{Martin-Lof1975},
and Prawitz \citeyearpar{Prawitz1977}. 
As these are well known,
the question arises
to what extent
their different ways of implementing
that general strategy
may suggest
alternatives for a Brouwerian
intuitionist to adapt and adopt.
I will look at
the theory of constructions
in Goodman’s version
as it was specifically adapted from Kreisel's
to avoid impredicativity.
The idea of canonical proofs
will be discussed in quite general terms
because,
for all the differences that exist
between the exact ways in which
Dummett,
Prawitz, 
and Martin-Löf
have developed it,
the fundamental differences with a Brouwerian
account arise at a common level.

Characteristic of all these accounts,
including Brouwer’s,
is that the notion of a totality of all possible proofs
is replaced by that of
a potentially infinite
collection of proofs
of different levels,
and an operator is introduced
that 
reduces
an arbitrary proof
to one
in that collection.
One of the differences between Kreisel’s original theory of constructions
and Goodman’s version is precisely that Kreisel entertained the thought of a reducibility operator
but did not use it \citep[pp.126–127, 2.215]{Kreisel1965}.
The reduced proofs are such
that
they can
be constructed 
‘from below’,
so that the explanation what 
a reduced proof of 
\(A\) 
is
will not be circular.
The explanations of implication
then avoid the notion of
an arbitrary proof of 
\(A\),
and instead appeal to that of an arbitrary reduced one.
What distinguishes these approaches from one another
is their underlying notion of proof and
the aspect of proofs that is the basis for the reduction procedure.

In the Brouwerian case,
the levels are the stages of the Creating Subject's
activity
and,
as on his conception proofs only exist
in the sense of demonstrations,
the reducibility operator is the trivial one.

In Goodman’s modification of
Kreisel’s
theory of constructions,
the levels are identified with 
‘maximal grasped domains’ 
\citep[pp.109–110]{Goodman1970}.
Goodman
does not elaborate much on what a ‘grasped domain’ is,
but he does say that it is a domain of constructions that
‘has been grasped as a totality’,
that the sense of ‘grasp' here is such that 
‘we cannot grasp the whole of the constructive universe, which is always
only a potential totality’,
and that a 
‘maximal’ 
grasped domain is one that
‘include[s] everything which is immediately understood when their elements are
understood’.
For example,
the natural numbers form,
by themselves,
a grasped domain,
but not a maximal one because
implicit in the species of natural numbers
is that of constructive numerical functions
\citep[p.110]{Goodman1970}.
There is a basic domain
\(B\),
which is level 
\(0\)
and includes the natural numbers;
and there is
an operation 
\(E\) 
to extend a level
so as to obtain the next:
\begin{quote}
Given any level 
\(L\),
we suppose that we can extend
\(L\)
to a new level containing all the objects of
\(L\),
all proofs about objects of 
\(L\),
and certain additional constructions […]
We emphasize that this is not a stratification by logical type,
but rather a stratification according to the subject matter of proofs.
\end{quote}
Clearly,
what does and does not belong to a maximal grasped domain is,
in its dependence on
what is 
‘immediately understood when the elements are understood’,
vague,
and prone to grow over time,
as with further experience 
we may come to find more
to be
‘immediately understood’.

The extension operation
\(E\) 
suffers from similar
problems.
Goodman is aware of this:
In the companion paper 
\citealt{Goodman1973a},
‘The arithmetic theory of constructions’,
he writes that
\begin{quote}
the rule which leads from the
\(n\)th 
level to the
\(n+1\)st
is not a rule which we can understand.
If it were,
then we could understand the notion of proof
in an absolute sense
and could visualize the entire constructive universe.
But that leads at once to self-reflexive paradoxes.
\end{quote}
(I will return to the role of visualisation below.)

In spite of this vagueness,
an operator
\(F\)
is introduced:
\begin{quote}
Suppose we have a grasped domain%
\footnote{[By this point,
Goodman has adopted the convention that
the grasped domains will always be understood to be maximal
\citep[p.110]{Goodman1970}.]}
\(a\).
Then we wish
\(a\)
to be as self-contained as possible.
Therefore,
if we have any rule
\(z\),
then we suppose that we can find a rule
\(Faz\)
in the domain
\(a\)
with the property that,
if
\(x\)
is in
\(a\)
and
\(zx\)
is defined and in
\(a\),
then
\[
\vdash Fazx \equiv zx.
\]
\end{quote}
It is crucial to the conception of levels as built from below,
Goodman points out,
that 
\(Fa\)
applied to
\(z\)
yields a rule already present in
\(a\)
that 
represents
\(z\)
in
\(a\)
\textit{extensionally};
\(z\)
may itself have been defined
at a much higher level,
but
the definition of
\(Faz\)
should not depend on that of
\(z\)
\citep[p.110]{Goodman1970}.

Now the key claim is this:
\begin{quote}
The introduction of the reducibility operator
\(F\)
is made necessary by the impredicative character
of intuitionistic implication.
It seems to us essential to the intuitionistic position
that given a fixed assertion 
\(\mathfrak{A}\)
about a well-defined domain,
there is always an a priori upper bound
to the complexity of possible proofs of
\(\mathfrak{A}\).
In case 
\(\mathfrak{A}\)
is an implication,
this principle already guarantees the existence of some sort of reducibility operator.
\citep[p.111]{Goodman1970}
\end{quote}
Specifically,
Goodman assumes that if a statement involves quantification over a domain
of level at most
\(p\),
then if there is a proof of that statement at all,
there is one of level
\(p+1\).
Under that assumption,
implication can be defined in a way that 
is not impredicative.
As mentioned,
such
a reducibility hypothesis had been 
stated 
by Kreisel, 
but not used 
(his theory of constructions is not stratified).
The hypothesis is a very strong one.
Being an existence claim,
it requires
that such a proof at level 
\(p+1\)
be obtainable
from
an arbitrary proof of
\(\mathfrak{A}\)
by a
constructive method.

However,
in
‘The arithmetic theory of constructions’,
Goodman states his conviction that
‘All of constructive mathematics is ultimately based
on finitary computation and on reasoning about
finitary computation’ 
\citep[p.284]{Goodman1973a}.
From what he goes on to say,
it is clear that for Goodman
‘reasoning about finitary computation’
proceeds by visualisation:
\begin{quote}
These insights can always be put into the following form:
a particular rule
\(a\),
applied to any element of this basic universe,
always gives the value
\(\underaccent{\mkern2.5mu\thicksim}{T}\).
Such an insight is evidently in
\(\forall\exists\)
form.
It asserts that a certain rule is actually a function~–
specifically,
the constant function whose value is always
$\underaccent{\mkern2.5mu\thicksim}{T}$.
We can think of the insight as the visualization,
or grasping,
of the totality of the computations of the values of the function.
The formal proof,
which is a finite object,
is not the insight but only a guide to aid us
in the visualization of this infinite structure of computations. 
[…]
Thus we can repeat the entire construction above
and consider insights which involve
visualizing the whole of level one.
These insights will be of level two.
\citep[pp.284,286]{Goodman1973a}
\end{quote}
Given this conviction
about what constructive mathematics consists in,
the existence of the required reducibility operator
may be easier to argue for,
and Goodman makes an attempt in that direction:
\begin{quote}
{}Let us 
suppose we have a rule
\(a\)
which is not of level 
\(0\).
Suppose we apply
\(a\)
to an object
\(b\)
of level
\(0\)
and obtain a result,
\(ab\),
which is again of level
\(0\).
Then it might happen that the computation of the value 
\(ab\)
involves the visualization of the whole of level
zero or of some even more complicated totality.
It seems to us,
however,
that that visualization cannot actually be essential to the computation.
For,
the visualization of such a totality of computations cannot make it possible
for us to make any computation other than the ones being visualized.
Moreover,
since the result of the computation does not involve
any infinite insight,
the computation of
\(ab\)
can only be using,
so to speak,
a bounded part of the infinite visualization.
In other words,
we are asserting that
if an infinite visualization is essential to a computation,
then it must occur in the result of the computation.
We may imagine,
for example,
that the computations are given in a normal form in which
the act of visualization plays the role of an
\(\omega\)-rule.
Then,
following Brouwer's proof of the bar theorem,
we are saying that any operation which plays the role of a cut-rule
can be eliminated,
so that only visualizations which occur in the result
are necessary for the computation.

It follows from this line of reasoning that if we are given a rule
\(a\)
which is not of level zero,
then we can think of it as a rule
\(c\)
of level zero by simply ignoring any request for a visualization
which cannot be carried out at level zero.
\citep[pp.285–286]{Goodman1973a} 
\end{quote}
Further clarification would be needed to explain
what,
if 
‘the visualization of a totality of computations cannot make it possible
for us to make any computation other than the ones being visualized’,
and if computation is all we are interested in,
would
\textit{motivate} 
introducing such a visualisation
into the construction of a rule
\(a\)
in the first place.
Be that as it may,
I take it that this,
or something like it,
is what Goodman had in mind
in the paper of 1970;
it would be hard to see in what other way
he could have hoped to justify the claim he makes.

As an understanding of constructive reasoning about proofs,
this limitation
to finitary reasoning
is overly restrictive,
as Kreisel
justly complained
in his review for 
\textit{Zentralblatt}
\citep{Kreisel1974}.%
\footnote{Kreisel will have thought of the introduction to 
\citealt{Godel1958}.}
A  more charitable reading
of the 1970 paper than Goodman allows
himself in his elucidation of 1973,
however,
is vulnerable to
the point made by Weinstein
that
\begin{quote}
The fact that various higher order intuitionistic systems
are conservative extensions of Heyting arithmetic with respect
to universal decidable sentences is rather weak evidence
for such a claim
[of the existence of the reducibility operator],
since it may just indicate our current inability to successfully
exploit the intuitionistic notions. 
\citep[p.266]{Weinstein1983}%
\footnote{Since these conservativity results had not been invoked by Goodman, 
it seems that they were the most plausible consideration that
Weinstein could think of
in favour of Goodman’s claim
(his list of references includes neither 
\citealt{Goodman1973a} 
nor  
\citealt{Kreisel1974}.)}
\end{quote}
In conclusion of this discussion of Goodman’s theory,
its notion of level
(maximal grasped domain)
is much less clear than in the Brouwerian account
(demonstrations effected up to a given stage of the Creating Subject’s activity),
and the required
reducibility operator is
highly problematic.

The other constructive account of implication
to consider is based on the notion of
‘canonical proof’,
introduced
by Dummett in his talk at the Bristol
Logic Colloquium in 1973
and published as 
‘The philosophical basis of intuitionistic logic’ 
\citep{Dummett1975}.
(Brouwer also had a notion of 
canonical proof,
but that is significantly different,
as will be discussed below.) 
A canonical proof is one in which the last step towards the conclusion
is the introduction of the main connective in the conclusion.
The key observation is that an inference such as

\begin{prooftree}
\AxiomC{$(A \rightarrow B) \rightarrow A$} 
\AxiomC{$A \rightarrow B$}
\RightLabel{$\rightarrow$-E}
\BinaryInfC{$A$}
\end{prooftree}

\noindent
and indeed
modus ponens
generally,
is not canonical.
Impredicativity of the 
clause of the Proof Interpretation
for implication
then would be avoided
by
explaining
a proposition
\(A \rightarrow B\)
not in terms of
\textit{any}
proof
of the antecedent,
but 
in terms of proofs 
that are canonical
\citep[p.394ff]{Dummett1977}.
Sundholm 
\citeyearpar[p.149]{Sundholm1994} 
has remarked that that 
requirement would be unnecessarily strong:
Implication may be defined for propositions whose proofs
\textit{can be 
reduced} 
to canonical form.
The function
\(f\)
that proves
\(A \rightarrow B\)
would then begin by performing that reduction.
Of course also this latter way of going about it
depends for its success on the existence of a proof
that is wholly canonical.

The reduction procedure would fulfill the role
of the reduction operator in the general scheme of
stratification outlined above;
the stratification would have only one level,
that of 
the canonical proofs.
The status of the reducibility thesis
that every proof is equivalent to one in canonical form
depends entirely on one’s exact notion of proof.
In Martin-Löf’s
type theory,
proofs are reducible to
canonical form
by definition:
in that theory,
that an object 
\(a\)
is a proof of proposition
\(A\)
\textit{means}
that
\(a\) 
is ultimately reducible to
a proof in which all steps are canonical.
The reduction of the former to the latter
is a mechanical procedure based on the syntax of the language.
Any object that does not allow for such canonisation is 
for that reason not a proof.
For those not accepting that approach,
the reducibility thesis is a very substantive one.
Dummett,
who never discussed Martin-Löf's work in detail,%
\footnote{Sundholm informed me that part of the motivation behind his 
‘Vestiges of realism’
\citeyearpar{Sundholm1994}
was to encourage Dummett to do so;
alas,
Dummett in his published reply did not seize this opportunity.}
let alone adopt particular ideas of it,%
\footnote{I will leave it to others to determine how easy or difficult this might be,
given Dummett's other theoretical commitments.
See,
for a start,
\citealt{Sundholm1994} 
and Dummett’s reply, 
\citealt{Dummett1994}.} 
never quite managed to find a justifiable alternative construal.
That the attempt, 
in the first edition of
\textit{Elements of Intuitionism},
to identify canonical proofs with
normal derivations in a formal system
was problematic
was acknowledged there by Dummett himself
\citeyearpar[pp.400–403]{Dummett1977},%
\footnote{See also 
\citealt[section 8]{Prawitz1987},
and
\citealt{Weiss1997}.}
and likewise he considered 
the reworked conception 
in 
\textit{The Logical Basis of Metaphysics}
of 1991 
‘very shaky’
already then
\citeyearpar[p.277]{Dummett1991b}.%
\footnote{On this point also 
\citealt[p.375]{Prawitz1994}
and 
\citealt[pp.205–209]{Sundholm1998}.}
In the second edition of 
\textit{Elements of Intuitionism}
\citeyearpar{Dummett2000b},
Dummett instead chose
to argue that impredicativity of 
implication is harmless.
The particular case he discusses is that of
\((B \rightarrow C) \rightarrow D\).
The ability to understand
proofs of the antecedent
would seem to require
an ability to survey of 
all possible proofs of
it,
because 
there is no
a priori limit to the complexity
of an operation transforming
proofs of 
\(B\)
into a proof of
\(C\).
This would threaten compositionality of meaning,
as understanding
\((B \rightarrow C) \rightarrow D\)
then presupposes
an understanding of
propositions that are
of the same or greater complexity
than its own.
But,
Dummett argues,
this is not so;
all that is really required
to understand
the antecedent
is
the ability to recognise 
of an operation that it is
effective and that it
will transform any proof of 
\(B\) 
into a proof of 
\(C\). 
‘We need not survey or circumscribe possible such operations
in any more particular way than this’ 
\citep[p.274]{Dummett2000b}.

Be that as it may,
there are very general reasons
why 
to a Brouwerian 
the meaning-theoretical approach to logic
exemplified in the work of Dummett, Prawitz, and Martin-Löf,%
\footnote{It should be noted that 
of these three,
it is in various ways Martin-Löf's approach that comes closest
to Brouwer's;
I have in mind here in particular the role of the difference between
proof objects and demonstrations,
the importance accorded to acts more generally,
and its functional character;
see also footnote~\ref{L003}.
But
I share Veldman's 
reservations having to do with the
homogeneity of the notion of type
and with the largely linguistic character
of the theory
\citep[pp.312–313]{Veldman1984};
I would say these are typical for
an ontological descriptivist faced with
meaning-theoretical foundations
\citep[section 5]{Sundholm.Atten2008}.}
is not attractive;
and the differences are of mathematical significance,
as seen in the example of the Bar Theorem,
of which there is a demonstration on a Brouwerian construal of that term
but not a meaning-theoretical one.%
\footnote{Constructivists of other stripes often accept the principle of bar induction
as primitive,
an alternative that Brouwer had indicated in 
\citealt[p.63n7]{Brouwer1927B}. 
As a curiosity,
I mention that Gödel once suggested that a normal form for all
mathematical proofs might be found,
and expressed the hope of proving the Bar Theorem from that
\citep[p.223]{Atten2015}.}
These general reasons have been discussed in
\citealt{Sundholm.Atten2008}.
We there characterised Brouwer's views
as ‘ontological descriptivist’ 
as opposed to meaning-theoretical;
an ontological descriptivist is someone
for whom ‘the correctness of a knowledge claim
is […] ultimately reduced to matters of ontology’
\citep[p.71]{Sundholm.Atten2008},
and meaning theory plays no ultimate role.
Here I should like to make,
in that light, 
some remarks on
Brouwer’s notion of canonicity and its differences
from the meaning-theoretical notion of canonicity.

Canonicity as appealed to by Brouwer
in 
‘\german{Definitionsbereiche}’
\citeyearpar[p.64]{Brouwer1927B},
is not that of a small
set of general steps into which any proof
can be analysed.
His notion is general with respect to structure only in the sense
that for him mathematical 
demonstrations
are species that are well-ordered;%
\footnote{More precisely,
they are ‘pseudo-well-ordered’:
a well-ordered structure of which the elements are subsequently
marked ‘empty’ or ‘full’,
depending on whether they play a role in the demonstration or not
(\citealt[p.64n8]{Brouwer1927B} 
and
\citealt[pp.46–47]{Brouwer1981A}).}
but they are not general with respect to what 
the canonical steps are.
They may vary from case to case,
as is clear from the example of the
the three canonical steps that figure in Brouwer’s analysis of demonstrations
that a tree is barred.
The only criterion for canonical steps  is that they be
‘immediately given to intuition’
\citep[p.64n8]{Brouwer1927B}%
\footnote{\german{‘der Intuition unmittelbar gegeben’}}.
Such intuitive givenness will depend on the structure of the specific objects
involved in the demonstration,%
\footnote{One may wish to generalise,
with Husserl,
the notion of syntax,
and distinguish between the
syntax of linguistic expressions and the syntax of
categorially formed
objects: 
\citealt[section 42b–d and \german{Beilage I}]{Husserl1974}; 
\citealt[section 11]{Husserl1976a}; 
\citealt[p.247n]{Husserl1985b}.}
and on the specific property to be established.
Here lies also the reason why
Brouwer's thesis
that
if there is a demonstration that a tree is barred,
then there is one in the canonical form he specifies,
does not presuppose a grasp of all possible means to establish the presence of a bar,
which is an open-ended collection.
The three canonical steps
Brouwer specifies are obtained not from such a grasp,
but derived directly from properties of the objects involved:
Any other way to show that all choice sequences passing through
a given node in the tree meet the bar
would imply the presence of restrictions on those sequences that are not
stated in the hypothesis of the Bar Theorem.
(This suggestion will be worked out in a later note.)

Brouwer's notion of canonicity then
is in particular not constrained by linguistic forms,
of which logical form would be an example.%
\footnote{\label{L003}%
As Kreisel pointed out,
for Brouwer the 
(introduction and elimination)
rules for the logical connectives
would not
\textit{define}
them,
but
\textit{describe}
(languageless)
mathematical operations
for which they have been recognised to be valid
\citep[p.122, 2.121]{Kreisel1965}.
In this particular respect,
functional interpretations of intuitionistic logic come closer
to what intuitionists mean than
meaning-theoretical accounts
depending on cut elimination
(such as Dummett's and Prawitz', but not Martin-Löf's):
From the intuitionistic point of view, 
the former track the original
languageless process of construction,
whereas the latter yield a linguistic object
that is no longer a direct representation
of that process.
(This footnote was inspired by
\citealt[p.231]{Kreisel.MacIntyre1982}.
I thank Göran Sundholm for his emphasis,
in discussion of this footnote,
that Martin-Löf's Constructive Type Theory is a functional interpretation.)}
This view is expressed already in Brouwer's dissertation:
\begin{quote}
While thus mathematics is independent of logic,
logic does depend upon mathematics;
in the first place intuitive logical reasoning 
is that special kind of mathematical reasoning
which remains if,
considering mathematical structures,
one restricts oneself to relations of whole and part;
the mathematical structures themselves are in no
respect especially elementary,
so they do not justify any priority of logical reasoning
over ordinary mathematical reasoning.
(\citealt[p.127]{Brouwer1907}, trl.~\citealt[p.73]{Brouwer1975})
\end{quote}
And even though Brouwer did not at that point have his
notion of canonical demonstration yet,
he was already at that stage
prepared to allow that
the structure of a mental demonstration
and
the structure of a linguistic representation of it
do not coincide:
\begin{quote}
Even in domains of mathematics where no relations of whole and part enter,
\textit{the relations which were in the mind are often transformed into relations of whole and part},
when they must be communicated verbally to other people;
hereby the usual language of mathematics is imbued with that of logical reasoning.
However,
this fact is due only to the centuries-old tradition of logical terms in language,
in connection with its limited vocabulary.
(\citealt[p.128n, emphasis Brouwer's]{Brouwer1907}, trl.~\citealt[p.73]{Brouwer1975})
\end{quote}
By the time of
the Bar Theorem,
Brouwer 
exploits
another structural difference:
demonstrations,
conceived of as mental objects,
are in general potentially infinite,
whereas expressions in the languages we actually use
are finite.
This contrast is brought up by Brouwer as his main argument,
in 1927,
against Hilbert's formalism
(quoted on p.~\pageref{L002}
above);
it likewise serves to draw a contrast to the 
meaning-theoretical tradition.%
\vspace{2\baselineskip}

\noindent
\textit{Acknowledgement.} 
I am grateful to
Walter Dean,
Hidenori Kurokawa,
Per Martin-Löf,
Joan Moschovakis,
Uday S.~Reddy,
Adrian Rezuş,
Göran Sundholm,
and Anne Troelstra
for comments and discussion.
Earlier versions
(without section 4)
were presented
at the
ASL meeting in Seattle in April 2017,
at the workshop
‘Critical Views of Logic’
in Oslo in August 2017,
at the 
‘\french{Journées SND}’
in Paris in September 2017
and at the meeting
‘Axiomatic Thinking’
in Lisbon
in October 2017.
I thank the organisers for the invitations,
and the audiences  
for their questions and comments.
Finally,
I thank the Institute for Advanced Study,
Princeton,
for permission to quote from
unpublished material in the
Gödel Papers.\vspace{2\baselineskip}

\bibliographystyle{rsl}
\bibliography{IntuitionisticImplication}

\end{document}